\documentclass[12pt]{article}   

\usepackage{amsmath}             
\usepackage{amsthm}
\usepackage{amsfonts}            
\usepackage{amssymb}             
\usepackage{graphicx}            

\usepackage{cite}                
\usepackage{geometry}            
\usepackage{booktabs} 

\usepackage{subfigure}
\usepackage{caption} 

\usepackage{algorithm}
\usepackage{algpseudocode}

\newtheorem{theorem}{Theorem}[section]
\newtheorem{definition}{Definition}[section]
\newtheorem{assumption}{Assumption}
\newtheorem{proposition}{Proposition}

\title{A novel viewpoint for Bayesian inversion based on the Poisson point process}
\author{Zhiliang Deng\textsuperscript{1 \thanks{dengzhl@uestc.edu.cn}},  
Zhiyuan Wang\textsuperscript{2}, 
Xiaomei Yang\textsuperscript{3 \thanks{yangxiaomath@swjtu.edu.cn}},
Xiaofei Guan\textsuperscript{2}
\\
\small \textsuperscript{1}University of Electronic Science and Technology of China, School of Mathematical Sciences \\
\small \textsuperscript{2} Tongji University, School of Mathematical Sciences \\
\small \textsuperscript{3}Southwest Jiaotong University, School of Mathematics}
\date{\today}                    

\begin{document}

\maketitle                       

\begin{abstract}

We present a novel Bayesian framework for inverse problems in which the posterior distribution is interpreted as the intensity measure of a Poisson point process (PPP). The posterior density is approximated using kernel density estimation, and the superposition property of PPPs is then exploited to enable efficient sampling from each kernel component. This methodology offers a new means of exploring the posterior distribution and facilitates the generation of independent and identically distributed samples, thereby enhancing the analysis of inverse problem solutions.

\end{abstract}

\section{Introduction}

Inverse problems, which involve estimating unknown parameters from observed data, are ubiquitous in scientific and engineering disciplines such as geophysics \cite{Sambridge2002, Zunino2023}, medical imaging \cite{Arridge1999, Arridge2009}. These problems are usually ill-posed, meaning that small errors in the observed data can lead to significant deviations in the estimated parameters. To address this challenge, the Bayesian approach has emerged as a powerful framework for solving inverse problems \cite{Stuart2010}. By incorporating prior knowledge and quantifying uncertainties through probability distributions, Bayesian methods provide a robust and interpretable solution to ill-posed inverse problems \cite{Dashti2017, Kekkonen2019, Latz2020, Latz2023, Stuart2010}.
The core idea of the Bayesian approach is to treat the unknown parameters as random variables and update their probability distributions based on observed data using Bayes' theorem. This results in a posterior distribution that encapsulates both prior beliefs and the information contained in the data. Unlike deterministic methods  \cite{Engl1996, Kirsch1996}, e.g., Tikhonov regularization and Landweber iteration, which yield a single point estimate, Bayesian methods provide a full probabilistic description of the solution, enabling the quantification of uncertainties and the assessment of confidence in the results \cite{Kaipio2005, Tarantola2005}.
One of the key strengths of the Bayesian framework is its flexibility. It allows the incorporation of diverse types of prior information, such as smoothness, sparsity, or physical constraints, through the choice of prior distributions \cite{Bui-Thanh2014, Calvetti2018, Kekkonen2019}. Additionally, Bayesian methods naturally handle noise in the data and provide a principled way to balance data fidelity with prior knowledge. This makes them particularly well-suited for complex and high-dimensional inverse problems where traditional methods may struggle.


In recent years, advancements in computational techniques, such as Markov chain Monte Carlo (MCMC) sampling and variational inference, have significantly enhanced the practicality of Bayesian methods for solving large-scale inverse problems \cite{Blei2017, Lu2017, Robert2004}. These developments have enabled the application of Bayesian approaches to a wide range of real-world problems, from reconstructing images in medical tomography to inferring subsurface structures in geophysical exploration \cite{Arridge2009, Fichtner2010, Iglesias2014}. 

In this paper, we introduce a novel perspective on Bayesian inverse problems. Specifically, we interpret the posterior distribution as a sampling distribution of a Poisson point process (PPP). By scaling it with the parameter of a Poisson random variable, we construct an intensity measure. The PPP associated with this intensity measure encapsulates information about the solution to the inverse problem. This approach allows us to explore the posterior distribution within the framework of a PPP, enabling the generation of posterior samples through realizations of the corresponding PPP.  Leveraging the specific structure of the Bayesian posterior distribution, we propose novel sampling algorithms for generating posterior samples. These algorithms are both simple and efficient.


The paper is organized as follows. In section 2 we make some general definitions about the PPP. In section 3 we discuss in detail the interpretation of Bayesian inverse problems in the PPP. 
In section 4 we provide some numerical tests to verify the proposed approach. In section 5 we present some conclusions.

\section{Poisson point process}
This section provides the foundational knowledge of PPPs used as preliminary background. The content covered here can be found in \cite{Last2017, Singh2009}, and additional details and perspectives are available in \cite{Daley2003, Jacobsen2006, Streit2010}.
Let \((\mathbb{X}, \mathcal{X})\) be a measurable space, and let \(\mathbf{N}_{<\infty}(\mathbb{X}) \equiv \mathbf{N}_{<\infty}\) denote the space of all measures \(\mu\) on \(\mathbb{X}\) such that \(\mu(B) \in \mathbb{N}_0=\mathbb{N}\cup \{0\}\) for all \(B \in \mathcal{X}\).  
A measure $\nu$ on $\mathbb{X}$ is said to be $s$-finite if $\nu$ is a countable sum of finite measures.
Let \(\mathbf{N}(\mathbb{X}) \equiv \mathbf{N}\) be the space of all measures that can be expressed as a countable sum of measures from \(\mathbf{N}_{<\infty}\).  
Let \(\mathcal{N}(\mathbb{X}) \equiv \mathcal{N}\) denote the \(\sigma\)-algebra generated by the collection of all subsets of \(\mathbf{N}\) of the form  
\[
\{\mu \in \mathbf{N}: \mu(B) = k\}, \quad B \in \mathcal{X}, \, k \in \mathbb{N}_0.
\]  
Let \((\Omega, \mathcal{F}, \mathbb{P})\) be an underlying probability space.

\begin{definition}
A point process  on \(\mathbb{X}\) is a random element \(\eta\) of \((\mathbf{N}, \mathcal{N})\), that is, a measurable mapping \(\eta: \Omega \to \mathbf{N}\).
\end{definition}

If \(\eta\) is a point process on \(\mathbb{X}\) and \(B \in \mathcal{X}\), then we denote by \(\eta(B)\) the mapping \(\omega \mapsto \eta(\omega, B) := \eta(\omega)(B)\). Obviously, we know that $\eta$ are random variables taking values in $\bar{\mathbb{N}}_0:=\mathbb{N}_0\cup \{\infty\}$, that is
\begin{align*}
\{\eta(B)=k\}\equiv\{\omega\in\Omega: \eta(\omega, B)=k\}\in\mathcal{F},\quad B\in\mathcal{X},\, k\in \bar{\mathbb{N}}_0.
\end{align*}
In this case we call $\eta(B)$ the number of points of $\eta$ in $B$. 

A point process \(\eta\), defined as a random measure on \((\mathbf{N}, \mathcal{N})\), inherently encodes the distribution of both the number of points and their spatial configuration. The total number of points, \(\eta(\mathbb{X})\), is a random variable whose distribution defines the probabilities \(\{p_n\}\), where \(p_n = \mathbb{P}(\eta(\mathbb{X}) = n)\). For each \(n\), conditioned on \(\eta(\mathbb{X}) = n\), the random measure \(\eta\) induces a symmetric probability measure \(\mathbb{P}_n\) on \(\mathbb{X}^n\). Specifically, for any measurable sets \(B_1, B_2, \dots, B_n \in \mathcal{X}\), \(\mathbb{P}_n\) assigns probabilities to the joint event that the \(n\) points lie in \(B_1, B_2, \dots, B_n\), respectively. This construction ensures that \(\mathbb{P}_n\) describes the spatial distribution of the points given their number. Thus, the structure of \(\eta\) naturally gives rise to the probabilities \(\{p_n\}\) for the number of points and the measures \(\{\mathbb{P}_n\}\) for their spatial configuration, providing a comprehensive description of the point process.
\begin{definition}
Let $\Lambda$ be an $s$-finite measure on $\mathbb{X}$. A PPP with intensity measure $\Lambda$ is a point process $\eta$ on $\mathbb{X}$ with the following two properties:
\begin{description}
\item[(i)] For every $B\in\mathcal{X}$, the distribution of $\eta(B)$ is Poisson with parameter $\Lambda(B)$, that is to say $\mathbb{P}(\eta(B)=k)=\frac{\Lambda(B)^k}{k!}e^{-\Lambda(B)}$ for all $k\in\mathbb{N}_0$; 
\item[(ii)] For every $m\in\mathbb{N}$ and all pairwise disjoint sets $B_1, \cdots, B_m\in\mathcal{X}$ the random variables $\eta(B_1), \cdots, \eta(B_m)$ are independent. 
\end{description}

\end{definition}
When the intensity measure $\Lambda$ admits a Radon-Nikodym derivate $\lambda$ with respect to a reference measure $\mu$ on $\mathbb{X}$, the function $\lambda$ is termed the intensity function of the point process.

\begin{definition}\label{def2.3}
Let $\nu$ and $\mu$ be probability measures on $\mathbb{N}_0$ and $\mathbb{X}$, respectively. Suppose that $\theta_1, \theta_2, \cdots$ are independent random elements in $\mathbb{X}$ with distribution $\mu$, and let $\kappa$ be a random variable 
with distribution $\nu$, independent of $(\theta_k)$. Then 
\begin{align}
\eta:=\sum_{k=1}^\kappa \delta_{\theta_k}
\end{align}
is called a mixed binomial process with mixing distribution $\nu$ and sampling distribution $\mu$. 
\end{definition}
The following result provides the key for the construction of Poisson processes.
\begin{proposition}\label{proposition}\cite{Last2017}
Let $\mu$ be a probability measure on $\mathbb{X}$ and let $\gamma\geq 0$. Suppose that $\eta$ is a mixed binomial process with mixing distribution $Poisson(\gamma)$ and sampling distribution $\mu$. Then $\eta$ is a PPP 
with intensity measure $\gamma\mu$.
\end{proposition}

When \(\mathbb{X}\) is a finite-dimensional space, (e.g. \(\mathbb{X} = \mathbb{R}^d\)),  the PPP admits an intuitive interpretation. To generate i.i.d. samples from a PPP with a non-trivially structured intensity function, 
  a two-step procedure can be employed (Algorithm~\ref{alg1}), see \cite{Streit2010} for details).  For low-dimensional spaces, Algorithm~\ref{alg1} provides an efficient sampling procedure. However, in general settings, significant computational challenges emerge
  when the intensity measure  $\Lambda(B)$, $B\in\mathcal{X}$, is analytically intractable, particularly when evaluating  $\Lambda(B)$ requires high-dimensional numerical integration.


The following theorems allow us to construct new PPPs from existing ones. The proofs can be found in \cite{Kingman1993, Last2017}. 
  \begin{theorem}[Superposition theorem] \label{superposition}
  Let $\eta_i$, $i\in\mathbb{N}$, be a sequence of independent PPPs on $\mathbb{X}$ with intensity $\Lambda_i$. Then
  \begin{align*}
  \eta:=\sum_{i}^\infty \eta_i
  \end{align*}
  is a Poisson process with intensity measure $\Lambda:=\Lambda_1+\Lambda_2+\cdots$. 
  \end{theorem}
  
    \begin{theorem}[Mapping theorem] 
  Let $\eta$ be a  PPP on $\mathbb{X}$ with intensity $\Lambda$ and let $\mathbb{Y}, \mathcal{Y}$ be a measurable space. Suppose $f:\mathbb{X}\rightarrow \mathbb{Y}$ be a measurable function. 
  Write $\Lambda^*$ for the induced measure on $\mathbb{Y}$ given by $\Lambda^*(B):=\Lambda(f^{-1}(B))$ for all $B\in\mathcal{Y}$. If $\Lambda^*$ has no atoms, then $f\circ \eta$ is PPP on $\mathbb{Y}$ with intensity measure $\Lambda^*$.
  \end{theorem}
  

For a proper PPP $\eta$ in Definition \eqref{def2.3}, one can assign to  each of the points $\theta_n$   a random mark $\varpi_n$ taking values in a measurable space $(\mathbb{M}, \mathcal{M})$, called the mark space. Let $\ell:\mathbb{X}\times \mathcal{M}\rightarrow [0, 1]$
be a probability kernel, i.e., $\ell(\theta, \cdot)$ is a probability measure for each $\theta\in\mathbb{X}$ and $\ell(\cdot, C)$ is measurable for each $C\in\mathcal{M}$. 
From the PPP \(\eta = \sum_{k=1}^\kappa \delta_{\theta_k}\), one can construct a marked point process 
\begin{align}\label{mark}
\varsigma:= \sum_{k=1}^\kappa \delta_{(\theta_k, \varpi_k)},
\end{align}
where, given \(\eta\), each \(\varpi_k\) is a random element in the mark space \((\mathbb{M}, \mathcal{M})\), independently drawn according to the conditional distribution \(\ell(\theta_k, \cdot)\). 
If the probability kernel takes the specific form  
\[
\ell_\pi(\theta, \cdot) := (1 - \pi(\theta)) \delta_0 + \pi(\theta) \delta_1, \quad \theta \in \mathbb{X},
\]  
where \(\pi: \mathbb{X} \to [0,1]\) is a measurable function, then the marked PPP in \eqref{mark}  induces a thinned PPP of \(\eta\), 
where each point \(\theta\) is independently retained with probability \(\pi(\theta)\) and removed with probability \(1 - \pi(\theta)\).

\begin{algorithm}
\begin{algorithmic}[1]
\State \textbf{Preliminaries} 
\begin{description}
\item[-] Let $\lambda_{\max}\geq \lambda(\theta)$, $\theta\in B\in\mathcal{X}$. 
\end{description}
\State \textbf{Step 1.} Give $n\in\mathbb{N}$. Generate candidate samples $\{\theta^{(1)}, \cdots, \theta^{(n)}\}$ according to a homogenous PPP with intensity $\lambda_{\max}$. 
\State \textbf{Step 2. Thinning process} 
For each $\theta^{(i)}$, compute the acceptance probability 
$$t=\frac{\lambda(\theta^{(i)})}{\lambda_{\max}}.$$
Retain the sample $\theta^{(i)}$ with probability $t$; otherwise, reject it. 
\end{algorithmic}
\caption{Realization of a PPP with thinning algorithm}
\label{alg1}
\end{algorithm}


\section{The perspective of PPP on Bayesian inversion}

%

Let \(\mathbb{X}\) and \(\mathbb{Y}\) be separable Banach spaces, each equipped with Borel \(\sigma\)-algebras \(\mathcal{X}\) and \(\mathcal{Y}\), respectively. Consider a measurable mapping
 \(G: \mathbb{X} \rightarrow \mathbb{Y}\), which typically originates from a system of differential equations.
We aim to find $\theta$ from a measurement of the form
\[
u = G(\theta) + \xi,
\]
where \(\xi \in \mathbb{Y}\) represents the noise. In this context, \(\theta\) refers to the parameter in the differential equation system, and \(u\) denotes the observation of the system's solution.
We assume that \((\theta, u)\) is a random variable in \(\mathbb{X} \times \mathbb{Y}\). Our objective is to determine the distribution of the conditional random variable \(\theta \mid u\) (i.e., \(\theta\) given \(u\)). This distribution allows us to estimate \(\theta\) and quantify the uncertainty associated with this estimation.

In the Bayesian framework, the random variable \(\theta|u\) is characterized by a prior measure \(\theta \sim \mu_{\rm pr}\) and a noise model \(\xi \sim \Xi\), where \(\mu_{\rm pr}\) (resp. \(\Xi\)) is a measure on \(\mathbb{X}\) (resp. \(\mathbb{Y}\)). 
The random variable $\theta|u$ is then distributed according to measure $\mu^u$ called the posterior measure. When $\mu^u\ll \mu_{\rm pr}$ for $\mu_{\rm pr}$-a.s., there exists a potential $\Phi:\mathbb{X}\times\mathbb{Y}\rightarrow\mathbb{R}$
\begin{align}\label{post}
\frac{d\mu^u}{d\mu_{\rm pr}}(\theta)=\frac{1}{Z}\exp(-\Phi(\theta; u)),
\end{align}
where $Z=\int_{\mathbb{X}} \exp\left(-\Phi(\theta)\right)d\mu_{\rm pr}(\theta)$.
When the mapping \(\Phi(\cdot; u) \in C(\mathbb{X}; \mathbb{R})\) for a fixed \(u\), it is a \(\mu_{\rm pr}\)-measurable function, and it satisfies the condition \(\mathbb{E}_{\mu_{\rm pr}}[\frac{1}{Z}\exp(-\Phi(\theta; u))] = 1\) \cite{Kekkonen2019}. The function \(\Phi(\cdot; u)\) is referred to as the negative log-likelihood, representing the fidelity to the data. For simplicity, we assume that $\Xi$ is a centered Gaussian measure with covariance operator $\Sigma$ on $\mathbb{Y}$. Consequently, the function $\Phi(\theta; u)$ takes the following form
\begin{align}
\Phi(\theta; u)=\left\|\Sigma^{-\frac{1}{2}}(G(\theta)-u)\right\|^2_{\mathbb{Y}}.
\end{align}
Moreover, we write the posterior measure as
\begin{align}\label{posteriormeasure}
d\mu^u(\theta)=\frac{1}{Z}\exp\left(-\Phi(\theta; u)\right)d\mu_{\rm pr}(\theta).
\end{align}
Obviously, $\mu^u$ is a probability measure on $\mathbb{X}$. 

In practical numerical implementations, the target parameter \(\theta\) is typically represented in a finite-dimensional basis, 
or the forward model is approximated using discretized algorithms such as finite element methods or spectral methods. This yields an approximated potential $\Phi$, denoted by $\tilde{\Phi}_N$. In particular, we define $\tilde{\mu}^N$ by
\begin{align}\label{apprpost}
\frac{d\tilde{\mu}^u_N}{d\mu_{\rm pr}}(\theta)=\frac{1}{\tilde{Z}^N}\exp\left(-\tilde{\Phi}^N(\theta; u)\right),
\end{align}
where $\tilde{Z}^N=\int_{\mathbb{X}} \exp\left(-\tilde{\Phi}^N(\theta; u)\right)d\mu_{\rm pr}(\theta)$. 
In \cite{Stuart2010}, A. Stuart established the convergence of the approximate measure \(\tilde{\mu}^u_N\)  to \(\mu^u\)  in Hellinger distance and total variation distance, under the following assumptions:
\begin{assumption}\label{assum1}
The function $\Phi:\mathbb{X}\times\mathbb{Y}\rightarrow\mathbb{R}$ has the following properties:
\begin{description}
\item[(i)] For every $\epsilon>0$ and $r>0$ there is an $M=M(\epsilon, r)\in\mathbb{R}$ such that, for all $\theta\in\mathbb{X}$ and all $u\in\mathbb{Y}$ with $\|u\|_{\mathbb{Y}}<r$,
\begin{align*}
\Phi(\theta; u)\geq M-\epsilon\|\theta\|_{\mathbb{X}}^2.
\end{align*} 
\item[(ii)] For every $r>0$ there is a $K=K(r)>0$ such that, for all $\theta\in\mathbb{X}$ and $u\in \mathbb{Y}$ with $\max\{\|\theta\|_{\mathbb{X}}, \|u\|_{\mathbb{Y}}\}<r$,
\begin{align*}
\Phi(\theta; u)\leq K.
\end{align*}
\end{description}

\end{assumption}

\begin{theorem}\label{thmStuart2010}\cite{Stuart2010}
Assume that the measures $\mu$ and $\tilde{\mu}^u_N$ are both absolutely continuous with respect to $\mu_{\rm pr}$, satisfying $\mu_{\rm pr}(\mathbb{X})=1$, with Radon-Nikodym derivatives given by \eqref{post} and \eqref{apprpost}
and that $\Phi$ and $\tilde{\Phi}^N$ satisfy Assumption \ref{assum1} with constants uniform in $N$. Assume also that for any $\epsilon>0$ there exists $K=K(\epsilon)>0$ such that
\begin{align}
|\Phi(\theta; u)-\tilde{\Phi}^N(\theta; u)|\leq K\exp(\epsilon\|\theta\|_{\mathbb{X}}^2)\psi(N),
\end{align}
where $\psi(N)\rightarrow 0$ as $N\rightarrow \infty$. Then the measures $\mu$ and $\tilde{\mu}^u_N$ are close with respect to the Hellinger distance (resp. the TV distance):
\begin{align}\label{conv}
\begin{aligned}
&d_{\rm Hell}(\mu^u, \tilde{\mu}^u_N):=\sqrt{\frac{1}{2}\int_{\mathbb{X}}\left(\sqrt{\frac{d\mu^u}{d\mu_{\rm pr}}}-\sqrt{\frac{d\tilde{\mu}^u_N}{d\mu_{\rm pr}}}\right)^2d\mu_{\rm pr}(\theta)}\leq C\psi(N),\\
&d_{\rm TV}(\mu^u, \tilde{\mu}^u_N):=\frac{1}{2}\int_{\mathbb{X}}\left|\frac{d\mu^u}{d\mu_{\rm pr}}-\frac{d\tilde{\mu}^u_{N}}{d\mu_{\rm pr}}\right|d\mu_{\rm pr}(\theta)\leq \sqrt{2}C\psi(N). 
\end{aligned}
\end{align}
\end{theorem} 

The posterior probability measure $\mu^u$ can be interpreted as defining the spatial distribution of a PPP. Building on Proposition \ref{proposition}, we:
\begin{description}
\item[i.]  Define a Poisson  measure $\nu$ with parameter $\gamma>0$
\item[ii.]  Construct a mixed binomial process $\eta$ by using:
\begin{description}
\item[-] $\nu$ as the mixing distribution;
\item[-] $\mu^u$ as the spatial sampling distribution.
\end{description}
\end{description}
This construction yields a PPP $\eta$ whose:
\begin{description}
\item[$\bullet$] Intensity measure is given by $\Lambda=\gamma\mu^u$, where $\gamma>0$ scales the expected number of points;
\item[$\bullet$] Spatial structure is entirely governed by $\mu^u$, with point locations sampled i.i.d. from $\mu^u$ conditional on the Poisson count. 
\end{description}
This probabilistic interpretation provides a powerful framework for analyzing the posterior measure $\mu^u$. 
Moreover, when we incorporate the approximate potential $\tilde{\Phi}_N$, we obtain an approximate PPP $\tilde{\eta}_N$ with corresponding intensity measure 
$\tilde{\Lambda}_N=\gamma\tilde{\mu}_N^u$. The finite-dimensional approximation $\tilde{\eta}_N$ maintains the PPP structure while enabling computationally tractable inference, with approximation error 
governed by the convergence properties of 
$\tilde{\mu}_N^u$ to the true posterior measure  $\mu_u$. 


Next, we discuss the convergence of $\tilde{\eta}^N$ to $\eta$ under the conditions of Theorem \ref{thmStuart2010}. First we have
\begin{align}\label{eqn8}
d_{\rm TV}(\tilde{\Lambda}_N, \Lambda)\leq \gamma d_{\rm TV}(\mu^u, \tilde{\mu}_N^u)\rightarrow 0 \quad \text{as}\quad N\rightarrow \infty. 
\end{align}

\begin{theorem}\label{thmvogue}
For $f\in C_b(\mathbb{X})$, if follows that
\begin{align}
\left|\int_{\mathbb{X}}f(\theta)d\tilde{\Lambda}_N(\theta)-\int_{\mathbb{X}}f(\theta)d\Lambda(\theta)\right|\leq \|f\|_\infty\cdot d_{\rm TV}(\tilde{\Lambda}_N, \Lambda) \rightarrow 0 \quad \text{as}\quad N\rightarrow \infty.
\end{align}
\end{theorem}
Theorem \ref{thmvogue} shows that the sequence of measures $\{\tilde{\Lambda}_N\}$ converges weakly to the measure $\Lambda$.  

\begin{theorem}
Let $A_1, \cdots, A_k\in\mathcal{X}$ be arbitrary pairwise disjoint compact Borel subsets of $\mathbb{X}$.  We have
\begin{align}\label{eqn10}
\left(\tilde{\eta}_N(A_1), \cdots, \tilde{\eta}_N(A_k)\right)\rightarrow_d \left(\eta(A_1), \cdots, \eta(A_k)\right),
\end{align}
where $\rightarrow_d$ indicates the convergence in distribution. 
\end{theorem}
\noindent {\bf Proof.} For all $i=1, \cdots, k$, $\tilde{\eta}_N(A_i)\sim Poisson(\tilde{\Lambda}_N(A_i))$ and $\eta(A_i)\sim Poisson(\Lambda(A_i))$. 
The convergence in \eqref{eqn8} implies that $\tilde{\Lambda}_N(A_i)\rightarrow \Lambda(A_i)$ as $N\rightarrow \infty$ for all $i=1, \cdots, k$. 
Moreover, for each $A_i$ we can obtain $\tilde{\eta}_N(A_i)\rightarrow_d\eta(A_i)$ as $N\rightarrow\infty$ by considering the convergence of the corresponding Laplace transform
\begin{align*}
\begin{aligned}
&\sum_{l=0}^\infty \frac{\tilde{\Lambda}_N(A_i)^l}{l!}e^{-\tilde{\Lambda}_N(A_i)}e^{-sl}=\exp(\tilde{\Lambda}_N(A_i)(e^{-s}-1))\\
&\rightarrow \exp(\Lambda(A_i)(e^{-s}-1))=\sum_{l=0}^\infty \frac{\Lambda(A_i)^l}{l!}e^{-\Lambda(A_i)}e^{-sl}, \quad \text{as} \quad N\rightarrow \infty.
\end{aligned}
\end{align*} 
Therefore, the joint characteristic function satisfies that
\begin{align}
\mathbb{E}\left[\exp\left(\sum_{i=1}^k \tilde{\eta}_N(A_i)(e^{-t_i}-1)\right)\right]\rightarrow\prod_{i=1}^k \exp\left(\Lambda(A_i)(e^{-t_i}-1)\right) \quad \text{as}\quad N\rightarrow \infty.
\end{align}
By L\'evy's continuity theorem, we get \eqref{eqn10}. 

\begin{theorem}
We have $\tilde{\eta}_N$ weakly converge to $\eta$, i.e., for any $f\in C_b(\mathbb{X})$, 
\begin{align}
\int_{\mathbb{X}}f(\theta)d\tilde{\eta}_N\rightarrow \int_{\mathbb{X}}f(\theta)d\eta \quad \text{as}\quad N\rightarrow \infty.
\end{align} 
\end{theorem}
\noindent {\bf Proof.}  According to Theorem \ref{thmvogue}, we get the result by considering the convergence of the corresponding characteristic function
\begin{align*}\begin{aligned}
&\mathbb{E}\left[e^{-\int_{\mathbb{X}}fd\tilde{\eta}^N}\right]=\exp\left(-\int_{\mathbb{X}}(1-e^{-f})d\tilde{\Lambda}_N\right)\\
&\rightarrow \exp\left(-\int_{\mathbb{X}}(1-e^{-f})d\Lambda\right)=\mathbb{E}\left[e^{-\int_{\mathbb{X}}f d\eta}\right] \quad \text{as} \quad N\rightarrow \infty, \quad \text{for all} \quad f\in C_b(X). 
\end{aligned}
\end{align*}


As established in \cite{Maroulas2020}, the intensity measure $\Lambda$ admits a density $\lambda$ with respect to a reference measure on the state space $\mathbb{X}$.  
Clearly, the reference measure can be chosen as the prior 
$\mu_{\rm pr}$, which implies that the intensity function takes the form
\begin{align}
\lambda=\frac{\gamma}{Z}\exp\left(-\Phi(\theta; u)\right).
\end{align}
It can be proven that $\eta$ is an $s$-finite measure. 
In certain cases—for example, when \(\mathbb{X} = \mathbb{R}^d\) and the prior measure \(\mu_{\rm pr}\) admits a density \(\pi_{\rm pr}\) with respect to the Lebesgue measure—the intensity function can be expressed as
\begin{align}  
\lambda(\theta) = \frac{\gamma}{Z} \exp\left(-\Phi(\theta; u)\right) \pi_{\rm pr}(\theta).
\end{align}  
In this case, $\eta$ is a $\sigma$-finite measure a.s..

\section{A decomposition sampling approach}



%

We now propose a novel sampling method for the posterior distribution by leveraging the superposition property of the PPP, as outlined in Theorem \ref{superposition}. In the subsequent analysis, we assume that the posterior (resp. prior) distribution defined in \eqref{posteriormeasure} admits a density $\pi^u(\theta)$ (resp. $\pi_{\rm pr}(\theta)$) with respect to the Lebesgue measure $d\theta$, i.e., $d\mu^u(\theta) = \pi^u(\theta)\,d\theta$ (resp. $d\mu_{\rm pr}=\pi_{\rm pr}d\theta$).

Our approach begins by approximating the posterior density using a Gaussian kernel mixture:
\begin{align}\label{mix_m}
\pi^u(\theta) \approx \sum_{k=1}^K w_k \varphi_k(\theta),
\end{align}
where each component $\varphi_k(\theta) = \mathcal{N}(\theta; \vartheta_k, \Xi_k)$ is a multivariate Gaussian density with mean $\vartheta_k$ and covariance matrix $\Xi_k$, and $\{w_k\}$ are non-negative weights satisfying $\sum_{k=1}^K w_k = 1$. 
 Multikernel learning algorithms have been extensively studied in the literature; see, for example, \cite{Gonen2011, Haugh2015}. In this work, we adopt a widely used approach—namely, the Expectation-Maximization (EM) algorithm (see Algorithm \ref{EM}) to obtain the approximation \eqref{mix_m}.

\begin{algorithm}
\begin{algorithmic}[1]
\State \textbf{Preliminaries} 
\begin{description}
\item[-] Draw $M$ samples $\{\theta^{(i)}\}$ from the prior distribution $\mu_{\rm pr}$;
\item[-] Compute $r^{(i)}=\Phi(\theta^{(i)})\pi_{\rm pr}(\theta^{(i)})$. 
\item[-] Normalization: $\tilde{r}^{(i)}=r^{(i)}/\sum_j r^{(j)}$.
\item[-] $w_k=1/K$ for all $k=1, 2, \cdots, K$. 
\end{description}
\State \textbf{EM:} 
\begin{description}
\item[-] E-step: Compute 
\begin{align}
\rho_{ik}=\frac{w_k\varphi_k(\theta^{(i)})}{\sum_l w_l \varphi_l(\theta^{(i)})}.
\end{align}
\item[-] Let $\hat{\rho}_{ik}=\tilde{r}^{(i)}\rho_{ik}$.
\item[-] M-step: Update
\begin{align}
w_k=\frac{\sum_i \hat{\rho}_{ik}}{\sum_{i, j}\hat{\rho}_{ij}}, \,\, \vartheta_k=\frac{\sum_i \hat{\rho}_{ik}\theta^{(i)}}{\sum_i\hat{\rho}_{ik}},\,\, \Xi_k=\frac{\sum_i \hat{\rho}_{ik}(\theta^{(i)}-\vartheta_k)(\theta^{(i)}-\vartheta_k)^\top}{\sum_i\hat{\rho}_{ik}}.
\end{align}
\end{description}
\State Repeat the EM step.
\end{algorithmic}
\caption{EM algorithm}
\label{EM}
\end{algorithm}
Each component in the approximation given by \eqref{mix_m} can be interpreted as the intensity function of a PPP. Using Algorithm~\ref{alg1}, we can independently sample from each component, thereby simplifying the overall sampling procedure.

\section{Numerical examples}
\subsection{Two-dimensional unimodal example}

We begin by considering the following one-dimensional elliptic boundary-value problem:
\begin{align*}
-\frac{d}{dx}\left(\exp(\theta_1)\frac{d}{dx}u(x)\right)=1,\,\, x\in[0, 1],
\end{align*}
with boundary condition $u(0)=0$ and $u(1)=\theta_2$. The forward problem involves computing the solution $u$ for given parameters $\theta_1$ and $\theta_2$, while the inverse problem aims to estimate  $\theta_1$ and $\theta_2$ from noisy measurements $u(\cdot; \theta_1, \theta_2)$ at locations $x_1=0.25$ and $x_2=0.75$. The explicit solution to this boundary-value problem is given by
\begin{align*}
u(x; \theta_1, \theta_2)=\theta_2 x+\exp(-\theta_1)\left(-\frac{x^2}{2}+\frac{x}{2}\right)
\end{align*}
and thus the forward map is defined as
\begin{align*}
G(\theta)=\left(u(x_1; \theta), u(x_2; \theta)\right)^\top\,\, \text{with}\,\, \theta=(\theta_1, \theta_2)^\top.
\end{align*}
The observation data are generated by corrupting the model output $G(\theta)$ with additive Gaussian measurement noise $\xi\sim\mathcal{N}(0, \Sigma)$, where the covariance matrix is $\Sigma=0.01\cdot I_2$ and $I_2\in\mathbb{R}^{2\times 2}$ denotes $2\times 2$ matrix. We adopt a Gaussian prior distribution $\mathcal{N}(0, I_2)$ for the parameters. This constitutes a nearly Gaussian example, originally introduced in \cite{Ernst2015} and subsequently employed in \cite{Garbuno-Inigo2020, Herty2019, Reich2021}. For our numerical experiments, we use the reference values $\theta^\dag=(0.0865, -0.8157)^\top$ and  $G(\theta^\dag)=(-0.0173, -0.573)^\top$ from \cite{Reich2021}. 

\begin{figure}
\centering
\subfigure[True pdf]{\includegraphics[width=0.4\textwidth]{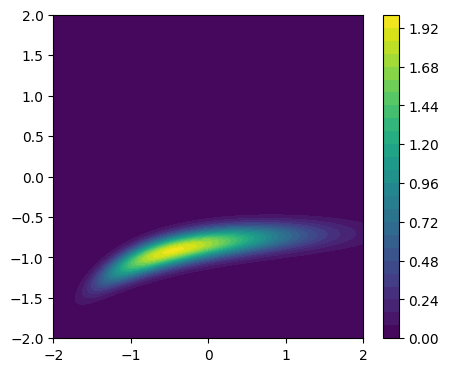}\label{ea}}
\subfigure[The estimated pdf and the samples]{\includegraphics[width=0.4\textwidth]{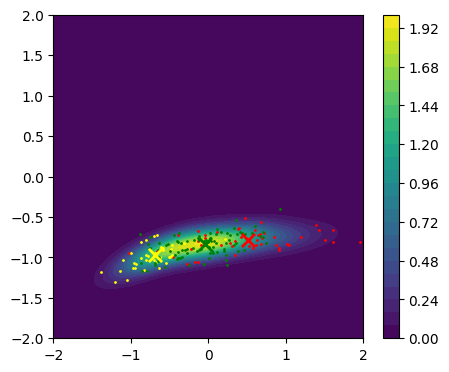}\label{eb}}
\caption{Numercial illustration for 2D unimodal example.}
\label{fig:ex1}
\end{figure}

In this example, we use $10,000$ prior samples to learn the posterior density with $3$ kernels, i.e., $K=3$ in \eqref{mix_m}.  
The true and estimated pdfs are shown in Figure \ref{fig:ex1}-\subref{ea} and \ref{fig:ex1}-\subref{eb}, respectively. Samples were generated from each component using the parameters listed in Table \ref{table_2d_uni} by Algorithm~\ref{alg1}. Figure \ref{fig:ex1} displays sample points colored by component, with their center positions from Table \ref{table_2d_uni} indicated by crosses. These color and cross markers are consistently used in the following examples.  These results show that the proposed method is effective.


\begin{table}
\centering
\begin{tabular}{cccc}
\hline
\textbf{Component} & $\boldsymbol{w}$ & $\boldsymbol{\vartheta}$ & $\boldsymbol{\Xi}$ \\
\hline
1 & 0.57507079 & $\begin{bmatrix} -0.01471933 \\ -0.83067489 \end{bmatrix}$ & $\begin{bmatrix} 0.19108289 & 0.03481514 \\ 0.03481514 & 0.02098995 \end{bmatrix}$ \\
\hline
2 & 0.25016956 & $\begin{bmatrix} -0.70383286 \\ -0.97366098 \end{bmatrix}$ & $\begin{bmatrix} 0.11865036 & 0.04646361 \\ 0.04646361 & 0.03302445 \end{bmatrix}$ \\
\hline
3 & 0.17475965 & $\begin{bmatrix} 0.71304739 \\ -0.77810972 \end{bmatrix}$ & $\begin{bmatrix} 0.27692851 & 0.0213613 \\ 0.0213613 & 0.01585085 \end{bmatrix}$ \\
\hline
\end{tabular}
\caption{Result parameters for 2D unimodal example.}
\label{table_2d_uni}
\end{table}

\subsection{Two-dimensional bimodal example.}
We next consider a two-dimensional bimodal example resulting from the nonlinear forward map
\begin{align}
G:\mathbb{R}^2\rightarrow\mathbb{R},\,\, G(\theta)=(\theta_1-\theta_2)^2. 
\end{align}
This example is given in \cite{Reich2021}. We use the same settings as in \cite{Reich2021}, i.e., the Gaussian prior with mean zero and covariance $I_2\in\mathbb{R}^{2\times 2}$ and Gaussian measurements errors $\xi\sim\mathcal{N}(0, I_2)$. The numerical results are based on the realization $G(\theta^\dag)+\xi=4.2297$ with $\theta^\dag=(-1.5621, -0.0021)^\top$.

In this example, we use $10,000$ prior samples to learn the posterior density with $2$ kernels, i.e., $K=2$ in \eqref{mix_m}. The result parameters are listed in Table~\ref{table2}. 
\begin{table}
\centering
\begin{tabular}{cccc}
\hline
\textbf{Component} & $\boldsymbol{w}$ & $\boldsymbol{\vartheta}$ & $\boldsymbol{\Xi}$ \\
\hline
1 & 0.498 & $\begin{bmatrix} 0.92807199 \\ -0.93907239 \end{bmatrix}$ & $\begin{bmatrix} 0.27007109 & 0.22820827 \\ 0.22820827 & 0.26950376 \end{bmatrix}$ \\
\hline
2 & 0.502 & $\begin{bmatrix} -0.93352663 \\ 0.93225046 \end{bmatrix}$ & $\begin{bmatrix} 0.27174892 & 0.22923468 \\ 0.22923468 & 0.26951865 \end{bmatrix}$ \\
\hline
\end{tabular}
\caption{Results parameters for 2D bimodal example.}
\label{table2}
\end{table}
\begin{figure}
\centering
\subfigure[True pdf]{\includegraphics[width=0.4\textwidth]{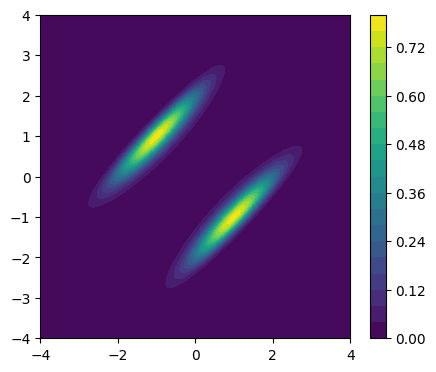}}
\subfigure[The estimated pdf and samples]{\includegraphics[width=0.4\textwidth]{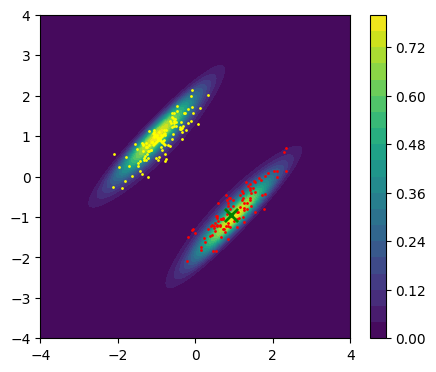}}
\caption{Numerical illustration for 2D bimodal example.}
\label{fig:ex2}
\end{figure}

\subsection{2D inverse heat conduction}
\begin{figure}
\centering
\subfigure[2D IHCP]{\includegraphics[width=0.3\textwidth]{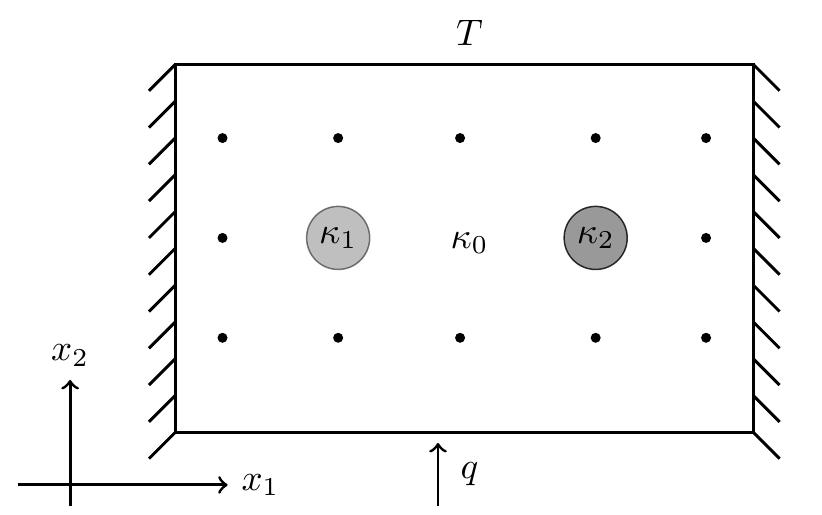}}
\caption{Heat conduction setups. }
\label{fig1}
\end{figure}

 We consider an example from \cite{Nagel2016}, which involves a stationary heat equation of the following form:
\begin{align}
\nabla\cdot(\kappa\nabla u)=0, \quad \text{in}\quad \Omega\subset \mathbb{R}^2, 
\end{align}
which describes a thermal distribution over a physical domain subject to appropriate boundary conditions \cite{Nagel2016}. 
We deal with the identification of thermal conductivities of inclusions in a composite material with close-to-surface measurements of the temperature.  The sketch setup of the  thermal problem is visualized in Fig. \ref{fig1}. 

The thermal conductivity of the 
background matrix is denoted as $\kappa_0$, while the conductivities of the material inclusions are termed as $\kappa_1$ and $\kappa_2$, respectively.  At the ``top" of the domain a Dirichlet boundary condition $u=T$ is imposed, while at the ``bottom" the Neumann
boundary condition $q=-\kappa_0\partial u/\partial y$ is imposed. Zero heat flux conditions $\partial u/\partial x=0$ are imposed at the ``left" and ``right" hand side. As the settings in \cite{Nagel2016}, the background thermal conductivity $\kappa_0$ 
is set to $15$,
while the true inclusion thermal conductivities are specified as $\kappa_1=32$ and $\kappa_2=28$. The boundary conditions $T=200$ and $q=2000$ are applied. 

\begin{figure}
\centering
\subfigure[The true pdf]{\includegraphics[width=0.4\textwidth]{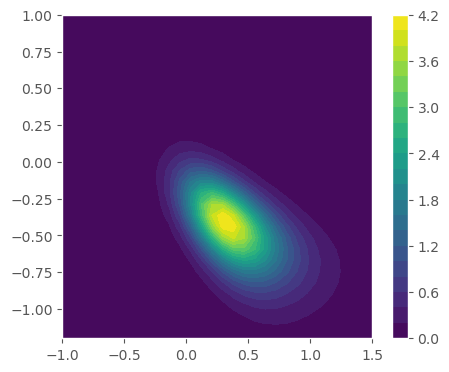}\label{4a}}
\subfigure[The estimated pdf and the Poisson samples]{\includegraphics[width=0.4\textwidth]{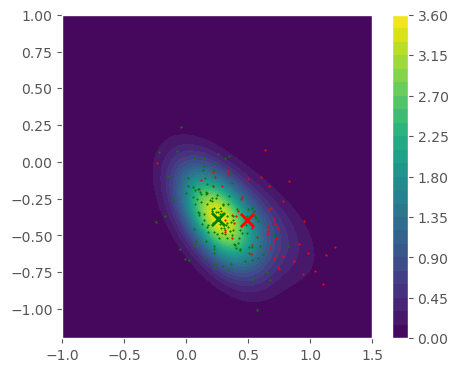}\label{4b}}
\caption{Numerical illustration for 2D inverse heat conduction. }
\label{fig2}
\end{figure}

In the reconstruction process, we transform \(\kappa_1\) and \(\kappa_2\) into new parameters using the relationship \(\kappa_{1, 2} = 30 + 6 \arctan \theta_{1, 2}\). Additionally, we assign a uniform prior distribution on the interval \([-1, 1]\) to \(\theta_{1, 2}\). The kernel density parameters are given in Table~\ref{2dh}. And the corresponding Poisson samples are collected in Figure~\ref{fig2}-\subref{4b}. From the displayed results, we can see that the Poisson sampling is suitable to this case.

\begin{table}
\centering
\begin{tabular}{cccc}
\hline
\textbf{Component} & $\boldsymbol{w}$ & $\boldsymbol{\vartheta}$ & $\boldsymbol{\Xi}$ \\
\hline
1 & 0.262 & $\begin{bmatrix} 0.49113893 \\ -0.39742986 \end{bmatrix}$ & $\begin{bmatrix} 0.08872554 & -0.03867398 \\ -0.03867398 & 0.05057162 \end{bmatrix}$ \\
\hline
2 & 0.738 & $\begin{bmatrix} 0.25307683 \\ -0.3884668 \end{bmatrix}$ & $\begin{bmatrix} 0.04205827 & -0.02528664 \\ -0.02528664 & 0.05140735 \end{bmatrix}$ \\
\hline
\end{tabular}
\caption{Result parameters for 2D IHCP.}
\label{2dh}
\end{table}

\subsection{High dimensional example}
\begin{figure}
\centering
\includegraphics[width=0.4\textwidth]{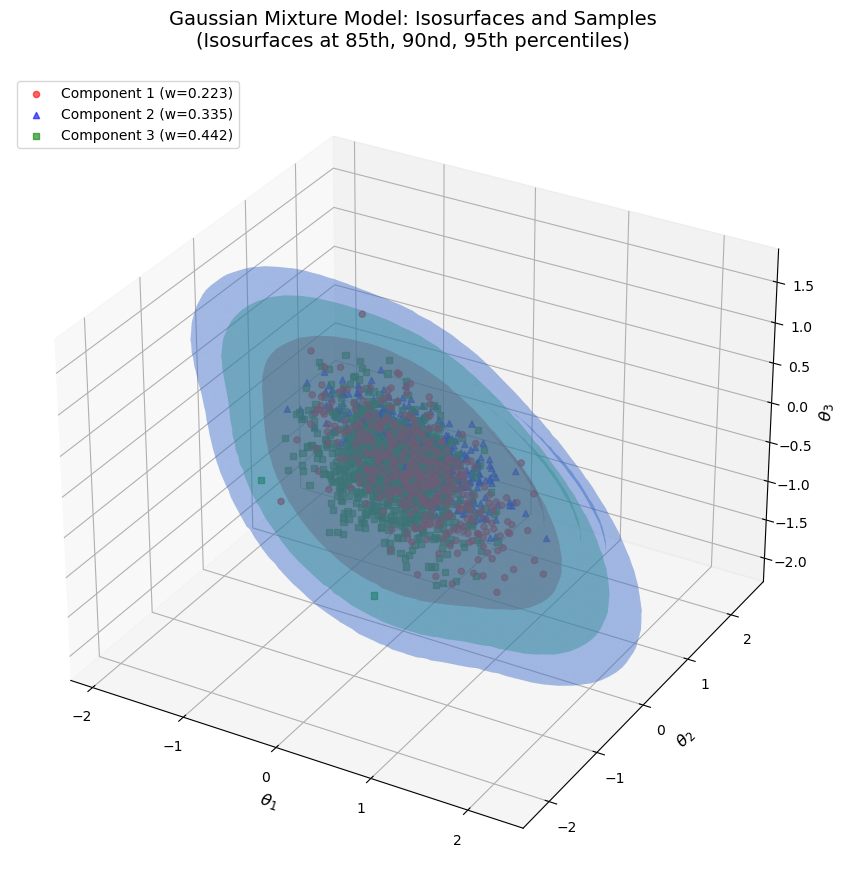}
\caption{Numerical illustration for high dimensional example. }
\label{fig5}
\end{figure}

We consider again the one-dimensional elliptic boundary-value problem
\begin{align}\left\{
\begin{aligned}
&-\frac{d}{dx}\left(\exp(p(x))\frac{d}{dx}u(x)\right)=1,&x\in[0, 1].\\
&u(0)=u(1)=0, &
\end{aligned}
\right.
\end{align}
The permeability function, $\exp(p(x))$, with $p\in L^\infty[0, 1]$, is to be determined from discrete observed data of $u$ within the interval $[0, 1]$. For any $p\in L^\infty[0, 1]$, there exists a solution $u\in H_0^1[0, 1]$. 
 We assume the prior for the unknown $p$ is Gaussian, $N(0, (-\frac{d^^2}{dx^2})^{-s})$, where the operator $-\frac{d^2}{dx^2}$ is defined with homogeneous Dirichlet boundary conditions. This prior admits a representation through the truncated Karhunen-Loève (KL) expansion:
 \begin{align}
 p(x) \approx \sum_{k=1}^N \lambda_k^{-s} \eta_k \phi_k(x):=\sum_{k=1}^N  \theta_k \phi_k(x)
 \end{align}
where $\{(\lambda_k, \phi_k)\}$ denotes the eigen-system of the aforementioned operator, and $\eta_k \sim N(0,1)$ are independent and identically distributed random variables. We take the true $\theta=[0.1, 0.4,-0.4]$ in this test.  We truncate the KL expansion to the first three terms for visualization purposes.
The GMM parameters  are given in Table~\ref{tab_3d}. The samples are illustrated in Figure~\ref{fig5}. We also plot the isosurfaces at 85th, 90nd, 95th percentiles in  Figure~\ref{fig5} respectivelly. It is obvious that the sample algorithm gives good results for this example. 
\begin{table}[htbp]
\centering
\begin{tabular}{cccc}
\toprule
\textbf{Component} & $\boldsymbol{w}$ &$\boldsymbol{\vartheta}$ &  $\boldsymbol{\Xi}$ \\
\midrule
1 & 0.223 & 
$\begin{bmatrix} 0.1505 \\ 0.0056 \\ -0.2299 \end{bmatrix}$ & 
$\begin{bmatrix}
 0.2408 & -0.0099 & -0.1583 \\
-0.0099 &  0.1432 &  0.0135 \\
-0.1583 &  0.0135 &  0.1284
\end{bmatrix}$ \\
\midrule
2 & 0.335 &
$\begin{bmatrix} -0.1205 \\ 0.4003 \\ -0.2407 \end{bmatrix}$ &
$\begin{bmatrix}
 0.1352 &  0.0218 & -0.0901 \\
 0.0218 &  0.1319 & -0.0429 \\
-0.0901 & -0.0429 &  0.0842
\end{bmatrix}$ \\
\midrule
3 & 0.442 &
$\begin{bmatrix} -0.0791 \\ -0.0685 \\ -0.2472 \end{bmatrix}$ &
$\begin{bmatrix}
 0.1360 & -0.0217 & -0.0848 \\
-0.0217 &  0.1738 &  0.0299 \\
-0.0848 &  0.0299 &  0.0703
\end{bmatrix}$ \\
\bottomrule
\end{tabular}
\caption{Result parameters  for high dimensional example.}
\label{tab_3d}
\end{table}


\section{Conclusion}
This paper has presented a novel approach to Bayesian sampling by formulating the posterior distribution within a PPP framework. The proposed method utilizes a thinning algorithm to draw samples from the PPP. By decomposing the intensity function into a Gaussian mixture and exploiting the superposition principle, efficient sample generation is achieved. Numerical results confirm the efficacy of the proposed method.

\bibliographystyle{plain}\bibliography{ref}

\end{document}